\input amstex
\documentstyle{amsppt}

\redefine\l{\lambda}

\pagewidth{10,5 cm}
\pageheight{15,5 cm}
\leftheadtext{}
\rightheadtext{}
\topmatter
\title
On a new inverse spectral problem
\endtitle
\author
A.M. Akhtyamov
\endauthor
\endtopmatter
\document
\subhead 1. Introduction \endsubhead
Inverse problems for the linear ordinary
differential operators containing parameter has been
studied by many authors
(\cite{1--7,10}). Operators pencils also
intensively studied (\cite{8--9}).
But this two directions has been developed independently.
The
differential operators
pencil recovery uniqueness theorem
will be proved
in our article. Novelty of this
result is not equation coefficient recovery,
but boundary conditions coefficients recovery.
We will show that all conditions of the theorem
are essential.
We will cite also a few interesting examples.

\subhead 2. The main result \endsubhead

Let us consider the two following
boundary--value problems
$$
\align
l(y,\l) &= y'' + b\l y' +c\l ^2y = 0,  \tag{1}\\
U_1(y)&=
a_{11}y(0)+
a_{12}y(1)+
a_{13}y'(0)+
a_{14}y'(1) = 0, \tag{2}\\
U_2(y)&=
a_{21}y(0)+
a_{22}y(1)+
a_{23}y'(0)+
a_{24}y'(1) = 0, \tag{3}
\endalign
$$
$$
\align
l(y,\l) &= y'' + b\l y' +c\l ^2y = 0, \tag{1} \\
 \widetilde U_1(y)&=
 \widetilde a_{11}y(0)+
 \widetilde a_{12}y(1)+
 \widetilde a_{13}y'(0)+
 \widetilde a_{14}y'(1) = 0, \tag{4}\\
 \widetilde U_2(y)&=
 \widetilde a_{21}y(0)+
 \widetilde a_{22}y(1)+
 \widetilde a_{23}y'(0)+
 \widetilde a_{24}y'(1) = 0, \tag{5}
\endalign
$$
Here

$\l $ is spectral parameter,

coefficients
$b,\ c,\
a_{ij},
\ \widetilde a_{ij}$ are complex numbers
and its independent on parameter~$\l$,

$rank(a_{ij})_{2 \times 4}=
rank(\widetilde a_{ij})_{2 \times 4}=2,$

$x \in [0,1],\ y = y(x) \in C^2[0,1].$

\proclaim{Theorem} If all nonzero eigenvalues of the
boundary--value problems
\ \thetag{1},
\thetag{2},
\thetag{3}
\ and
\ \thetag{1}, \thetag{4}, \thetag{5}
coincide, their multiplicities coincide, and,
in addition,
the following conditions
are realized :
\roster
\item $ b^2-4c\ne 0 $,
\item $ b\ne 0 $,
\item $ c\ne 0 $,
\endroster
\noindent
then spectral problems themselves coincide,
that is, the linear formes
$ \widetilde U_1(y),\
 \widetilde U_2(y) $ are linear expressed
with the help of the linear formes
$ U_1(y),\
 U_2(y).$
\endproclaim

\demo{Proof}

Let:

$g(\l)$ be an arbitrary entire function,

$l$ be an arbitrary integer number,

$y_1(x,\l )=e^{\omega _1 \l x},\
 y_2(x,\l )
=e^{\omega _2 \l x}$ be a fundamental solves system of the
equation
\thetag{1},

$\omega _1,\ \omega _2 $ be
roots of the characteristic equation $\omega ^2 + b\omega + c = 0$
(according to condition of the theorem
$\omega _1 \ne \omega _2 $)

$\widetilde\Delta(\l )$ be characteristic
equation~\thetag{1},\,\thetag{4},\,\thetag{5}.

\smallpagebreak
Under the condition
$\l \ne 0$ the characteristic determinant of
the following boundary--value problem
$$
\align
l(y,\l) &= 0, \tag{1} \\
 \widetilde U_1(y)&=
  0, \tag{4}\\
\l ^le^{g(\l)}(
 \widetilde U_2(y))&=
 0, \tag{6}
\endalign
$$
is the following function
$$
\align
\widetilde\Delta_1(\l )
&\equiv
\vmatrix
\widetilde U_1(y_1(x,\l )) &
U_1(y_2(x,\l )) \\
\l ^le^{g(\l)}
\widetilde U_2(y_1(x,\l )) &
\l ^le^{g(\l)}
\widetilde U_2(y_2(x,\l ))
\endvmatrix \\
&\equiv
\l ^le^{g(\l)}
\widetilde\Delta(\l ).
\tag{7}
\endalign
$$ \par

 It follows from \thetag{7} that  nonzero eigenvalues
of spectral problems
\thetag{1},\,\thetag{4},\,\thetag{5} and
\thetag{1},\,\thetag{4},\,\thetag{6}
coincide. Then according to the condition of the
theorem,
the nonzero eigenvalues of the spectral problems
\thetag{1},\,\thetag{2},\,\thetag{3}
and \thetag{1},\,\thetag{4},\,\thetag{6} also coincide.

Nonzero eigenvalues of the problem
\thetag{1},\,\thetag{2},\,\thetag{3}
are the roots of the following entire function
$$
\Delta(\l )
\equiv
\vmatrix
U_1(y_1(x,\l )) &
U_1(y_2(x,\l )) \\
U_2(y_1(x,\l )) &
U_2(y_2(x,\l ))
\endvmatrix
$$
(\cite{6, p.27}).

It follows from Weierstrass theorem
about an entire function representation by its roots that
$$
\widetilde\Delta_1(\l )
\equiv
\l ^k e^{f(\l )}
\Delta(\l ), \tag{8}
$$
where $f(\l)$ is a certain entire function,
and $k$ is a certain integer number.
We assume, that
$g(\l ) \equiv f(\l ), \ l=k$
at the equation \thetag{6}.
Then it follows from \thetag{7} and \thetag{8}
that
$$\Delta(\l ) \equiv\widetilde\Delta(\l ).$$
That is
$$
\vmatrix
U_1(y_1(x,\l )) &
U_1(y_2(x,\l )) \\
U_2(y_1(x,\l )) &
U_2(y_2(x,\l ))
\endvmatrix
\equiv
\vmatrix
\widetilde U_1(y_1(x,\l )) &
\widetilde U_1(y_2(x,\l )) \\
\widetilde U_2(y_1(x,\l )) &
\widetilde U_2(y_2(x,\l ))
\endvmatrix .
$$
From here
$$
\allowdisplaybreaks
\align
&(a_{11}a_{22}-a_{21}a_{12}-
\widetilde a_{11}\widetilde a_{22}+\widetilde a_{21}\widetilde a_{12})
(y_1(0)y_2(1)-y_1(1)y_2(0))\\
+&(a_{11}a_{23}-a_{13}a_{21}-
\widetilde a_{11}\widetilde a_{23}+\widetilde a_{13}\widetilde a_{21})
(y_1(0)y_2'(0)-y_1'(0)y_2(0))\\
+&(a_{11}a_{24}-a_{21}a_{14}-
\widetilde a_{11}\widetilde a_{21}+\widetilde a_{21}\widetilde a_{14})
(y_1(0)y_2'(1)-y_1'(1)y_2(0))\\
+&(a_{23}a_{12}-a_{13}a_{22}-
\widetilde a_{23}\widetilde a_{12}+\widetilde a_{13}\widetilde a_{22})
(y_1(1)y_2'(0)-y_1'(0)y_2(1))\\
+&(a_{12}a_{24}-a_{22}a_{14}-
\widetilde a_{12}\widetilde a_{24}+\widetilde a_{22}\widetilde a_{14})
(y_1(1)y_2'(1)-y_1'(1)y_2(1))\\
+&(a_{13}a_{24}-a_{23}a_{14}-
\widetilde a_{13}\widetilde a_{24}+\widetilde a_{23}\widetilde a_{14})
(y_1'(0)y_2'(1)-y_1'(1)y_2'(0))
\equiv
0.
\endalign
$$

If we substitute the functions
$$y_1(x,\l )=e^{\omega _1 \l x},\
 y_2(x,\l )
=e^{\omega _2 \l x}$$ for
the preceding identity, we obtain
$$
\allowdisplaybreaks
\align
\left(
\vmatrix
a_{11} & a_{12} \\
a_{21} & a_{22}
\endvmatrix
\right.
&-
\left.
\vmatrix
\widetilde a_{11} & \widetilde a_{12} \\
\widetilde a_{21} & \widetilde a_{22}
\endvmatrix
\right)
(e^{\omega _2 \l } -
e^{\omega _1 \l })
\\
+ \left(
 \vmatrix
a_{11} & a_{13} \\
a_{21} & a_{23}
\endvmatrix
\right.
&-
\left.
\vmatrix
\widetilde a_{11} & \widetilde a_{13} \\
\widetilde a_{21} & \widetilde a_{23}
\endvmatrix
\right)
(\omega _2 \l  -
\omega _1 \l )
\\
+\left(
\vmatrix
a_{11} & a_{14} \\
a_{21} & a_{24}
\endvmatrix
\right.
&-
\left.
\vmatrix
\widetilde a_{11} & \widetilde a_{14} \\
\widetilde a_{21} & \widetilde a_{24}
\endvmatrix
\right)
(\omega _2\l e^{\omega _2 \l } -
\omega _1\l e^{\omega _1 \l })
\\
+\left(
\vmatrix
a_{12} & a_{13} \\
a_{22} & a_{23}
\endvmatrix
\right.
&-
\left.
\vmatrix
\widetilde a_{12} & \widetilde a_{13} \\
\widetilde a_{22} & \widetilde a_{23}
\endvmatrix
\right)
(\omega _2\l e^{\omega _1 \l } -
\omega _1\l e^{\omega _2 \l })
\\
+\left(
\vmatrix
a_{12} & a_{14} \\
a_{22} & a_{24}
\endvmatrix
\right.
&-
\left.
\vmatrix
\widetilde a_{12} & \widetilde a_{14} \\
\widetilde a_{22} & \widetilde a_{24}
\endvmatrix
\right)
(\omega _2\l e^{(\omega _1+\omega _2) \l } -
\omega _1\l e^{(\omega _1+\omega _2) \l })
\\
+\left(
\vmatrix
a_{13} & a_{14} \\
a_{23} & a_{24}
\endvmatrix
\right.
&-
\left.
\vmatrix
\widetilde a_{13} & \widetilde a_{14} \\
\widetilde a_{23} & \widetilde a_{24}
\endvmatrix
\right)
\omega _1 \omega _2{\l}^2 (e^{\omega _2 \l } -
e^{\omega _1 \l })
\equiv
0.
\endalign
$$

By the hypotheses of the Theorem, we have
$$\omega _1 +\omega _2 \ne 0,\ \
\omega _1 -\omega _2 \ne 0,\ \
\omega _1 \ne 0,\ \ \omega _2 \ne 0. $$
Therefore functions
$$
\align
e^{\omega _2 \l } -
e^{\omega _1 \l }&,\ \ \
\omega _2 \l  -
\omega _1 \l,\\
\omega _2\l e^{\omega _2 \l } -
\omega _1\l e^{\omega _1 \l }&,\ \ \
(\omega _2\l e^{\omega _1 \l } -
\omega _1\l e^{\omega _2 \l }),\\
\omega _2\l e^{(\omega _1+\omega _2) \l } -
\omega _1\l e^{(\omega _1+\omega _2) \l }&,\ \ \
\omega _1 \omega _2{\l}^2 (e^{\omega _2 \l } -
e^{\omega _1 \l })
\endalign
$$
are linear independent functions with
respect to argument $\l $.
(It is easily verified by the definition of functions linear
independence.)
Conequently $$\allowdisplaybreaks\align
\vmatrix
a_{11} & a_{12} \\
a_{21} & a_{22}
\endvmatrix
&=
\vmatrix
\widetilde a_{11} & \widetilde a_{12} \\
\widetilde a_{21} & \widetilde a_{22}
\endvmatrix , \ \ \
\vmatrix
a_{11} & a_{13} \\
a_{21} & a_{23}
\endvmatrix
=
\vmatrix
\widetilde a_{11} & \widetilde a_{13} \\
\widetilde a_{21} & \widetilde a_{23}
\endvmatrix , \\
\vmatrix
a_{11} & a_{14} \\
a_{21} & a_{24}
\endvmatrix
&=
\vmatrix
\widetilde a_{11} & \widetilde a_{14} \\
\widetilde a_{21} & \widetilde a_{24}
\endvmatrix , \ \ \
\vmatrix
a_{12} & a_{13} \\
a_{22} & a_{23}
\endvmatrix
=
\vmatrix
\widetilde a_{12} & \widetilde a_{13} \\
\widetilde a_{22} & \widetilde a_{23}
\endvmatrix , \\
\vmatrix
a_{12} & a_{14} \\
a_{22} & a_{24}
\endvmatrix
&=
\vmatrix
\widetilde a_{12} & \widetilde a_{14} \\
\widetilde a_{22} & \widetilde a_{24}
\endvmatrix , \ \ \
\vmatrix
a_{13} & a_{14} \\
a_{23} & a_{24}
\endvmatrix
=
\vmatrix
\widetilde a_{13} & \widetilde a_{14} \\
\widetilde a_{23} & \widetilde a_{24}
\endvmatrix .
\endalign
$$
It follows from preceding equations that
all third order minors of the matrix
$$
A =
\pmatrix
a_{11} & a_{12} &
a_{13} & a_{14} \\
a_{23} & a_{24} &
a_{23} & a_{24} \\
\widetilde a_{13} & \widetilde a_{14} &
\widetilde a_{13} & \widetilde a_{14} \\
\widetilde a_{23} & \widetilde a_{24} &
\widetilde a_{23} & \widetilde a_{24}
\endpmatrix
$$
is equal to zero.
That is why rank of the matrix $A$ is equal to two.
The last means that the forms
$ \widetilde U_1(y),\
 \widetilde U_2(y) $ are linear expressed by forms
$ U_1(y),\
 U_2(y).$  As was to be proved.
\enddemo

\remark{Remark}
One need two, three or more spectrums of the especially
choosed problems for the differential operators recovery
uniquness. (\cite{2--6}). But according to our theorem
one needs only one spectrum for
the operators pencil recovery uniquness.
It shows that "pencil case" differs from
"nonpencil case".
\endremark

\subhead 3. On essence of each condition of the Theorem \endsubhead

We will show now that all conditions of the Theorem
are essential. We cite three examples, showing
the importance of each condition of the theorem.

\example{Example 1
({\rm First condition of the Theorem is missed:\ }
$b^2-4c=0,\ b \ne 0,\ c \ne 0$)}

The spectral problems
$$
\xalignat 2
&y''-2 \l y' + \l ^2 y = 0,  &&y''-2 \l y' + \l ^2 y = 0, \\
&y(0)=0,  &&y(0)=0, \\
&y'(1)=0,  &&y(1)+2y'(1)=0
\endxalignat
$$
have the same set of the eigenvalues.
This set consists of only one eigenvalue
$\l = -1. $
However
forms of boundary conditions of the
the first spectral problem are not
linear expressed by the boundary conditions
forms of the second spectral problem.
\endexample

\example{Example 2
({\rm Second condition of the Theorem is missed:\ }
$b^2-4c \ne 0,\ b = 0,\ c \ne 0$)}

The spectral problems
$$
\xalignat 2
&y''- \l ^2 y = 0, &&y''- \l ^2 y = 0, \\
&y(0)+2y'(0)=0, &&y(0)=0, \\
&y(1)=0, &&y(1)-2y'(1)=0
\endxalignat
$$
have the same set of the eigenvalues,
which is the same as roots of characteristic
determinant
$(1+2\l )e^{-\l} +
(-1+2\l )e^{\l}. $
However
forms of boundary conditions of the
the first spectral problem are not
linear expressed by the boundary conditions
forms of the second spectral problem.
\endexample

\example{Example 3
({\rm Third condition of the Theorem is missed:\ }
$b^2-4c \ne 0,\ b \ne 0,\ c = 0$)}

The both spectral problems
$$
\xalignat 2
&y''- \l ^2 y = 0, &&y''- \l ^2 y = 0, \\
&y(0)+2y'(0)=0, &&y(0)=0, \\
&y(1)=0, &&y(1)-2y'(1)=0
\endxalignat
$$
both have not eigenvalues. Therefore the sets of
the eigenvalues of the problems are coincided.
However
forms of boundary conditions of the
the first spectral problem are not
linear expressed by the boundary conditions
forms of the second spectral problem.
\endexample

\subhead 4. About some
generalisations of the Theorem
\endsubhead

If coefficients of the boundary conditions
depend on parameter $\l$, then conclusion of
the Theorem generally speaking is not true.
We cite an example,
confirming this preposition.

\example{Example 4
({\rm Coefficients of boundary conditions depend on parameter $\l$})}

The spectral problems
$$
\xalignat 2
&y''- 3 \l y' + \l ^2 y = 0, &&y''- 3 \l y' + \l ^2 y = 0, \\
&\l y(0)+y'(0)=0, &&\l y(0)+2y'(0)=0, \\
&4 \l y(1)+ y'(1)=0, &&2\l y(1)+y'(1)=0
\endxalignat
$$
have the same eigenvalues, which is the same as the roots
of the characteristic determinant
$\l ^2 (12e^{2 \l} -
15e^{\l}). $
However forms of
boundary conditions
of
the first spectral problem are not
linear expressed by the boundary conditions
forms of the second spectral problem.
If coefficients of boundary conditions
depend on parameter
$\l$, then they can not be univalent recovered by the spectrum.
It means that a generalisation of the Theorem on this way is not
possible.
\endexample

One can not univalent recover both coefficients of a equation
and coefficients of all boundary conditions.
We cite a corresponding example.

\example{Example 5
({\rm Coefficients of equations are not coincided})}

The spectral problems
$$
\xalignat 2
&y''+ 2 \l y' + \l ^2 y = 0, &&y''+ 4 \l y' + 5 \l ^2 y = 0, \\
&y(0)=0, &&y(0)=0, \\
&y(\pi )=0, &&y(\pi )=0
\endxalignat
$$
have the same eigenvalues
$\ \{\l _n\}=\pi n,\ \ n \in \Bbb Z .$
Boundary conditions of
the first spectral problem coincide with the boundary conditions
of the second spectral problem.
But coefficients of the both equations are different.
Thus one can not univalent recover both the coefficients
of the equation and the coefficients of the boundary
conditions.
\endexample

\enddocument
\end